\documentstyle[12pt]{article}
\def\={\!=\!}

\def\={\!=\!}

\def\|{{\Vert}}

\def\noi{\noindent}

\def\skip{\vspace{2mm}}
\newcount\u
\def\Remark{\skip\noi{{\bf{Remark \number\u.}}} \advance\u by 1}

   \title{The Common Knowledge of Formula  Exclusion}
\author {Robert Samuel Simon}
\begin{document}
\maketitle

\thispagestyle{empty}

\vfill

\noi London School of Economics\newline 
Department of Mathematics\newline 
Houghton Street\newline 
 London WC2A 2AE\newline 
R.S.Simon@lse.ac.uk\newline 
+44-207-955-6753 \newline fax +44-207-955-6877
       
\vfill     
   
\setcounter{page}{-1}
\noi           This research  was supported by 
     the Center for Rationality and Interactive Decision Theory (Jerusalem),  
    the Department of Mathematics of the Hebrew University (Jerusalem),    
 the Edmund Landau Center for Research in Mathematical
Analysis (Jerusalem)  
       sponsored by the Minerva Foundation (Germany), and  the German Science 
 Foundation (Deutsche Forschungsgemeinschaft).

\newpage
 \noi 
{\bf Abstract:} A  
 Kripke structure for the S5 logic  
 is defined by three   sets $S$, $J$ and $X$,  a collection
 $ ({\cal P}^j\ | \ j\in J)$ 
 of partitions of $S$ and a function $\psi : S \rightarrow \{ 0,1\}
 ^X$.   To each partition ${\cal P}^j$ 
 corresponds a person $j\in J$ who 
 cannot distinguish between any two points belonging to
 the same  member of ${\cal P}^j$ 
 but 
 can distinguish between different members of ${\cal P}^j$. A
 cell  is a minimal subset $C$ of $S$ such
that  for all $j$ the property 
 $P\in {\cal P}^j$ and $P \cap C \not= \emptyset$ implies  that 
 $P\subseteq C$.   Construct a sequence ${\cal R}_0, {\cal
 R}_1, \dots $ of partitions of $S$ inductively 
 by ${\cal R}_0 = \{ \psi^{-1} (a)
 \ | \ a\in \{ 0,1\} ^X\}$ 
 and $x$ and $y$ belong  to the same member of ${\cal
 R}_i $ if and only if $x$ and $y$ belong to the same member of ${\cal
 R}_{i-1}$ and for every person $j$ the members $P_x$ and $P_y$ 
 of ${\cal P}^j$
 containing $x$ and $y$ respectively  intersect the same members of
 ${\cal R}_{i-1}$.  Let ${\cal R}_{\infty}$ be the limit of the ${\cal
 R}_i$, namely $x$ and $y$ belong to the same member of ${\cal
 R}_{\infty}$ if and only if $x$ and $y$ belong to the same member of 
 ${\cal R}_i$ for every $i$. 
 For any sets $X$ and $J$  of persons there  is a  
canonical Kripke structure  defined on a set 
 $\Omega= \Omega (X,J)$ 
    such
 that   from any Kripke structure for the S5 logic
  using the same $X$ and $J$  there is a canonical map   to
 $\Omega$ with the property 
 that $x$ and $y$ are mapped to the same point of $\Omega$ if and
 only if $x$ and $y$ share the same member of ${\cal R}_{\infty}$.  
We define a cell of $\Omega$    to be  surjective if every  
 Kripke structure for the S5 logic   
  that maps to it does so surjectively. A cell of a Kripke structure 
 for the S5 logic 
 has finite fanout if every $P\in {\cal P}^j$ contained in the cell 
 is a finite set. 
   All cells of $\Omega$ with 
 finite fanout are surjective, but the converse does not hold.

 \vskip2cm 
 \noi 

 \thispagestyle{empty}

    \newpage

\section{Introduction}   
{\em Common knowledge} by all persons in $J$  
 of the  event $E$ means that  
for  every string of
persons  $i_1, i_2, \dots , i_k$ in $J$ it holds that $i_k$ knows that
$i_{k-1}$ knows that ...  $i_1$ knows that the event $E$ has occurred
(Lewis 1969).  
 One way to formalise   knowledge and  common
knowledge
 is   through 
semantic  models called Kripke structures,
  (see also Aumann, 1976.) In this paper, we assume throughout  
 the S5 logic (defined below), so that we will refer to
 these structures   simply as  
 Kripke structures. 
    Applying the above definition of common knowledge to this context (as 
 defined in the abstract), 
 a subset $A$ is known in common by all the
  persons in $J$  at the point $x\in A$  if the
 cell containing $x$ is contained in the set $A$. 
\vskip.2cm

\noi First, lets look at an example that
 illustrates part of the 
  problem solved in this paper. Let ${\bf \overline N}$ be 
 $\{ 1,2, \dots \} \cup \{ \infty\}$, and give 
 ${\bf \overline N}$ the topology where all integers are isolated 
 points however  the sequence $1$, $2$, ... converges to $\infty$.  Let 
 $S$ be ${\bf \overline N} \times {\bf \overline N}$ with the corresponding 
 product topology. Let $J$ be $ \{ 1,2,3\}$. For every $i$ define the set     
   $P_i^1= \{ (i,1), (i,2)
 , \dots , (i,\infty)\}$, including the possibility of $i=\infty$. 
  For every $i$ defined the   set 
 $P^2_i= \{ (1,i), (2,i), \dots ,  (\infty, i)\}$, including the
 possibility  of $i=\infty$.   For every $i$ define  the set $P^3_i =
 \{ (k,l) \ | \ k+l= i\}$ with $ P_{\infty}^3 = 
\{ (k,l) \ | \ $ either $k$ or $l$ is equal 
 to $\infty \}$.  
For each $j=1,2,3$ define the partition  ${\cal P}^j$ of $S$ to be 
 $\{ P_1^j, P_2^j , \dots , P_{\infty}^j\}$. Let $X$ be 
 a singleton $\{ x \}$. Let the evaluation function $\psi$  give the    
 evaluation 
  $1$ to the point  $(1,1)$, the evaluation  $1$ to the point $(2,1)$, 
  and everywhere else the 
evaluation  $0$.  It is not difficult to see that 
 the ${\cal R}_i$ as defined in the abstract  eventually 
 separates all the points of $S$ (once the three points $(0,0)$, $(1,0)$ and 
 $(0,1)$ are distinquished from each other).   According to results 
 presented later in  
 this paper, this Kripke structure  is homeomorphic to 
 a cell of $\Omega (X,J)$ that does not have finite fanout. 
 As we will see later, this cell  is  also  not surjective. Experimentation 
 with such examples gives the false impression that surjectivity
 and finite fanout are equivalent properties.  We will 
 discover later    why it is difficult to construct  such a  cell 
 that  is both  surjective and fails to have finite fanout.  We return to 
 this example  later.
  \vskip.2cm 
\noi The above example shares much in common with well known examples of 
 game theory.
In particular,  the two-person non-zero-sum  Electronic Mail Game of 
 Rubinstein (1989) 
 is based on a  structure of information  with  similar properties to 
 the above example. In the Rubinstein game,  
 there is a special point at infinity 
 where the players would  have common knowledge of  
 a certain  payoff matrix (the payoff matrix could be  determined by 
 an  evaluation function $\phi$) and it is  the 
 limit point of an infinite sequence of points that are  included in the game.
 However 
 this special 
  point at infinity  is excluded from the game and at the initial 
   point of the sequence 
  and only at this initial point    
 a different payoff matrix is valid. The ``almost common knowledge'' in 
 the title of the article refers to the almost common knowledge 
 of the payoff matrix valid at both  
  the special point at infinity and all but one of the other points.  
 The analysis of the Rubinstein  
game shows that the  equilibrium behavior at these limiting  
 points is very different from that at the excluded limit point at infinity.  
  The structure 
 of the Rubinstein example 
  has finite fanout, but by adding an extra third player to this structure 
 Fagin, Halpern and Vardi (1991) constructed an interesting 
 example of a cell of 
 the appropriate $\Omega$ 
that is not surjective and does not have finite fanout. Mapping to this cell 
 injectively are  
  two alternative Kripke structures, one where the third person 
 cannot distinguish between the special point and all the other points 
  and another 
 where the third person does exclude this special point.  As we will see 
 later, the point $(\infty,\infty)$ of our example is a  similar  special 
 point at infinity.\vskip.2cm 
\noi   The Electronic Mail Game shows that equilibrium behavior 
 can fail to be continuous with respect to changes 
in the information structure.
 A more radical but related discovery is a  three-person 
  non-zero-sum Bayesian game such that all  equilibrium 
 behavior  is not measurable with respect to the structure of information
 (Simon 2003).  
 Although the underlying 
  information structure of this Bayesian game is related closely  
 to the main result of this paper, 
  it   does have finite 
 fanout and therefore does not provide an example of a
 surjective cell lacking 
 finite fanout. 
    \vskip.2cm

\noi 
Before we can describe our main result, we must define $\Omega (X,J)$,
the canonical Kripke structure  . 
   \vskip.2cm 

\noi Let $X$ be a  set of primitive propositions, and let $J$ be a 
  set of agents. Although it is legitimate to consider the case of 
either  $X$ or  $J$ infinite, for this paper we will assume
  throughout that both $X$ and $J$ are finite.  
     Construct the set ${\cal L} (X,J)$ of  formulas using  the  
     sets $X$   
     and   $J$ in the
      following  way:\newline
      1) If $x\in X$ then $x\in {\cal L}(X,J)$,\newline
      2) If $g\in {\cal L}(X,J)$ then $(\neg g)\in {\cal L}(X,J)$,\newline
      3) If $g,h\in {\cal L}(X,J)$ then $(g\wedge h) \in 
      {\cal L}(X,J)$,\newline
      4) If $g\in {\cal L}(X,J)$ then $k_j g\in {\cal L}(X,J)$
for every
       $ j\in J$, \newline 
       5) Only formulas constructed through application of the 
        above four  rules are members of ${\cal L}(X,J)$.\newline 
  We write simply ${\cal L}$ if there is no ambiguity. $\neg f$ stands for
 the negation of $f$, $f\wedge g$ stands for  both  $f$ and $g$. 
 $f \vee g$ stands for either $f$ or $g$  (inclusive) and
 $f\rightarrow g$ stands for $\neg f \wedge g$.   
\vskip.2cm 

\noi If ${\cal K} =  (S, {\cal P}, J, X, \psi)$ is a Kripke structure
    then 
  define a map  $ \alpha^{{\cal K}}$ from ${\cal L}(X,J)$ to 
$ 2^S$, the subsets of $S$, 
 inductively on the structure of the formulas:
      \newline {\bf Case 1 $f=x\in X$:}  $ \alpha^{{\cal K}} (x):=  \{
  s\in S\ | \ \psi^x  (s) =1\}$.
     \newline {\bf Case 2 $f=\neg g$:} $ \alpha^{{\cal K}} (f) := 
     S\backslash \alpha^{{\cal K}}(g)$,
     \newline {\bf Case 3 $f=g\wedge h$:} $ \alpha^{{\cal K}} (f) := 
      \alpha ^{{\cal K}} (g) \cap  \alpha^{{\cal K}}
     (h)$,
     \newline {\bf Case 4 $f=k_j(g)$:} $ \alpha^{{\cal K}}(f):= 
     \{ s\ |\ s\in P\in {\cal P}^j 
\Rightarrow P\subseteq \alpha^{{\cal K}} (g)\}$.
\vskip.2cm 

\noi There is a very elementary  logic defined on the formulas in ${\cal
         L}$  called $S5$. 
For a longer  discussion of the $S5$ logic, see   
 Cresswell and Hughes (1968); 
 and for the multi-person  variation, 
  see Halpern and Moses (1992)       
  and also Bacharach, et al, (1997).
      Briefly, the $S5$ logic  
    is defined by two rules of inference, modus ponens and necessitation,  
    and  five types of axioms. 
      Modus ponens means that if $f$ is a theorem and 
       $f\rightarrow g$ is a theorem, then $g$ is also a theorem.   
       Necessitation 
      means that if $f$ is a theorem then $k_j f$ is also a theorem 
       for all $ j\in J$. The axioms are the following, for every   
        $f,g\in {\cal L}(X,J)$ and $ j\in J$:  
   \newline 1) all formulas resulting from theorems of the propositional  
    calculus through substitution, 
   \newline 2) $(k_j f\wedge k_j (f\rightarrow g) ) \rightarrow k_j g$,
   \newline 
     3) $k_j f \rightarrow f$, 
     \newline 4) $k_j f\rightarrow k_j (k_j f)$, 
     \newline 5) $\neg k_j f\rightarrow k_j (\neg k_j f)$. 
  \vskip.2cm 
  
\noi 
    A set of formulas  ${\cal A}\subseteq 
    {\cal L} (X, J)$  is called  {\em complete} 
    if for every
       formula $f\in {\cal L} (X,J)$  either $f\in {\cal A}$ 
       or $\neg f\in {\cal A}$.
      A set of
      formulas is called {\em consistent}
      if no finite subset
 of this set leads to a logical contradiction, meaning 
  a deduction of $f$ and  $\neg f$ for some formula $f$.
We define $$ \Omega (X,J):=\{ S\subseteq {\cal
L}(X,J)\ |\ 
 S\mbox{ is complete and consistent} \}.$$  
   
\noi 
$\Omega (X,J)$ is itself a  Kripke structure with evaluation. 
      For every person  $j\in J$ we define its corresponding  partition  
      ${\cal Q}^j(X,J)$ to be that 
  generated by the inverse images of  
 the function
  $\beta ^j: \Omega(X,J) \rightarrow 2^{{\cal L} (X,J)}$  namely 
 $$\beta ^j (z):= \{ f\in {\cal L} (X,J)\ | \ 
  k_j f \in z\}, $$ the set of formulas known by person $j$.     
  Due to the fifth set of axioms $\beta ^j (z)\subseteq \beta ^j (z')$ 
   implies that $\beta ^j (z) =\beta ^j (z')$. 
  We will write $\Omega$, ${\cal L}$ and ${\cal Q}^j$ if there is 
   no ambiguity.

     \vskip.2cm 
\noi   If ${\cal K} = (S, {\cal P}, J, X, \psi)$  is the corresponding 
 Kripke structure we  
      define a map $\phi^{{\cal K}}:S\rightarrow \Omega (X,J)$ 
       by
      $$\phi^{{\cal K}}(s):= \{ f \in {\cal L}(X,J)\ | \  \ 
     s\in  \alpha^{{\cal K}} (f)\}.$$ 
This is the canonical map referred to in the abstract, also contained
     in Fagin, Halpern, and Vardi (1991). 
\vskip.2cm 

\noi As stated in the abstract, a  cell $C$ of a
 Kripke structure   
 has {\em finite fanout}  if every choice of $i\in J$ and 
 $P\in {\cal P}^i$  contained in $C$ the set $P$  has finitely many elements. 
In Simon (1999) a cell of $\Omega$  
was defined to be 
 {\em surjective} if all Kripke structures     ${\cal K}$  
 that map to it by $\phi^{\cal K}$  do so surjectively. 
We construct an example of a countable and 
  surjective cell of  $\Omega$
  that 
  does not have finite fanout.
 (In Simon (1999) it was proven that any cell of $\Omega$ with finite
 fanout is surjective and   any surjective  cell of
 $\Omega$   must be 
 countable.)

\vskip.2cm 
\noi 
Central to understanding the relation between surjectivity and 
 finite fanout  is  
 point-set topology. For every Kripke structure   
  ${\cal K} =  (S, {\cal P}, J, X, \psi)$
  we define a topology
 on the set $S$, the same as in Samet (1990).  
 Let $\{  \alpha^{\cal K}(f)\ |\ f\in {\cal L}\} $ 
 be the base of open sets of $S$. We call this the topology {\em induced 
 by the formulas}.
 The topology of a subset $A$ of $S$ will be the relative  
  topology for which the open sets of $A$ are $\{ A\cap O\ |\ O$ is an open 
   set of $S\} $.    
      \vskip.2cm

\noi 
Why is our main result surprising?  It is closely related to
 representations of Kripke structures 
 through canonical structures
 indexed by  ordinal numbers.
\vskip.2cm 

  \noi 
 Fagin (1994) defined  for any  ordinal number
 $\gamma$ (and  a sets of persons  and 
 primitive propositions)   a hierarchically constructed  canonical 
  Kripke structure   
  $W_{\gamma}$ such that $W_{\omega}$ is $\Omega$,
(where $\omega$ stands for the first infinite ordinal).
 This canonical structure represented 
 all the possible 
 truth evaluations with the  ordinal 
 numbers  
  representing   the levels in the construction of these statements.   
There are alternative   canonical  constructions   
 corresponding to the ordinal numbers (Heifetz and Samet 1998, 1999), but with 
 regard to the first infinite ordinal $\omega$ they are the same as Fagin's.   
 For every  Kripke structure  
  and ordinal number $\gamma$ 
   there are  canonical 
maps  defined 
 to the canonical structures $W_{\gamma}$.    
\vskip.2cm 
  
\noi 
 If there is an ordinal $\alpha$ such that the map to 
  $W_{\alpha}$  is injective, then 
 the Kripke structure     is called {\em non-flabby}, and the 
  first such ordinal  is called the {\em distinguishing}
 ordinal. Otherwise the distinguishing ordinal is the first ordinal
 $\alpha$ where all pairs of points  which get mapped eventually to different
 places do so when mapped  to  $W_{\alpha}$.   
 There is another 
 minimal ordinal ${\beta}$, possibly larger than the distinguishing ordinal,
 for which the  image of the Kripke structure   
 in $W_{\beta}$ 
can be extended to any $W_{\gamma}$ with $\gamma > \beta$ 
 in only  one  unique way.  This  ordinal is called the  
  {\em uniqueness} ordinal. Fagin (1994) proved that the 
 uniqueness ordinal is a limit ordinal
 and  never greater than the next limit ordinal 
 above the distinguishing ordinal.
   Fagin established that the     
 necessary and sufficient condition  for a cell of $\Omega$   
  to have the first infinite ordinal $\omega$ as its 
 uniqueness ordinal is that the cell has 
     finite fanout.     
Without explicitly mentioning  topology, Fagin (1994) showed that any 
 member  $P$ of some ${\cal Q}^j$  
 is a compact subset of  $\Omega$.      
An extension to $W_{\omega +1}$ of an $z$ in $\Omega$ 
 is defined by  dense
 subsets $R_j$ of  the various  $P_j \in {\cal Q}^j$ containing $z$. Therefore 
  there is a unique extension of a cell of $\Omega$   if and only if 
  for every person $j$  
every $P_j\in {\cal P}^j$ in the cell has only one dense subset, 
 which is equivalent to the cell having finite fanout.  
  \vskip.2cm 
\noi If a cell of $\Omega$ does not have finite fanout, we know that 
 there is no unique extension of this cell to the higher levels. It is 
 plausible to believe that this lack of a unique extension can be realised 
 by alternative Kripke structures   that map injectively to $\Omega$. For the 
 structure  that maps injectively but not surjectively to this cell,    
   the persons would   have   
 common knowledge  that some set of formulas valid somewhere 
  in the cell are not valid at any point 
 in the original Kripke structure. 
  The surjective property is exactly the impossibility  
  of such   a 
  common knowledge of formula exclusion. 
\vskip.2cm

 \noi There is a  good reason why one can believe easily 
 that  the  
  surjectivity and finite fanout properties of cells are equivalent. 
 The relationship between the properties rests largely 
  upon a property 
  called  {\em centeredness}. 
 The centered property has several equivalent 
 definitions; the most straightforward definition is that a cell of 
 $\Omega$ is centered if and only if no other cell of 
$\Omega$ 
  shares the same  set of formulas 
 held in common knowledge (Simon 1999).    
  (The set of formulas held in common knowledge 
 is a constant throughout any given cell; see Halpern and Moses 1992). 
An equivalent formulation  of centeredness 
 is that the cell is an open set relative 
 to the closure of itself.  The 
 difference between centered and uncentered cells is radical; if a cell is 
 not centered  then there are uncountably many other cells 
 sharing the same set of formulas in common knowledge (Simon 1999).  
Furthermore  the converse does hold for centered cells of $\Omega$,
namely that 
  a  centered cell of
 $\Omega$ is surjective if and only if it has finite fanout (Theorem 3b, 
Simon  1999).  \vskip.2cm 

\noi
 The  lack of finite fanout for a cell $C$ of $\Omega$ 
implies the existence a cluster point $y$ 
 of some   $P\in {\cal Q}^j$ that is contained in $C$.  
Is the point $y$ is a  good candidate for  
  the existence of a Kripke structure    that maps to  
 $C\backslash \{ y\}$?  If $C$ is centered there will be such a
 Kripke structure   
that maps to $C$ but avoids the point $y$. 
\vskip.2cm 

\noi Returning to the above example, we see that the 
corresponding  
 cell of \newline  $\Omega (\{ x, \} , \{ 1,2,3\} )$ is  
  centered, a context in which surjectivity and finite fanout are equivalent. 
 According to Simon (1999) there is a  Kripke structure   that  
  maps injectively but not surjectively to this cell which 
 results from   removing the  
  single point  $(\infty, \infty)$.  According to the same theory, the 
 same can be done by  
 removing all the  cluster points, namely the set ${\cal P}_{\infty}^3$, and 
 furthermore these are the only two ways to map injectively but 
 not surjectively to this cell.   
 \vskip.2cm

\noi 
In the next section, we provide some more  background  necessary 
 to understand our solution.  In the third and concluding section, we 
 present our example of  a cell of $\Omega$ that is
 surjective but without finite fanout.

\section {More Background}

           Central to this paper is the first part of 
{\bf Lemma 5} of Simon (1999), which 
 states that  if ${\cal K}=(S;J; ({\cal P}^j\ |\ j \in J); 
    X;\psi)$ is a  Kripke structure   and    
    $P$ is a member of $ {\cal P}^j$ for some $j\in J$ then  
    $\phi^{\cal K}(P)$ 
     is a dense subset of $F$ for some $F\in {\cal Q}^j$. 
 This fact was used implicitly by Fagin (1994).   
\vskip.2cm 
   
\noi 
Given a 
 Kripke structure   ${\cal K}=(S; J ;
({\cal P}^j\ |\ j \in J);X; \psi )$ and a subset $A\subseteq S$,  
          we define the   
           Kripke structure   ${\cal V}^{{\cal K}}(A):= 
 (A;J ; 
({\cal P}^j|_A  \ |\ j \in J);X; \psi|_A )$ 
 where for all $j\in J$ 
 ${\cal P}^j|_A:= \{ F\cap A\ |\ F\cap A\not= \emptyset $ and $F\in  
  {\cal P}^j\} $. 
  We define a subset $A\subseteq \Omega$ to be  {\em good}  
       if 
       for every $j\in J$ and every 
       $F\in {\cal Q}^j$ satisfying $F\cap A\not= 
        \emptyset$  it follows that  $F\cap A$ is dense in $F$. 
        By Lemma 6 of Simon (1999)  $A$ is good if and only if  
     for every $z\in A\ $ $\phi^{{\cal  V}^{{\cal K}} (A)}(z) =z$.  
\vskip.2cm

      \noi  The  next lemmatta 
 relate directly the good property 
 to our problem.
\vskip.2cm

\noi  {\bf Lemma 7}  of Simon (1999): If 
   ${\cal K} = (S; J; ({\cal P} ^j | j\in J); X;\psi)$ is a   
    Kripke structure     
    then $  \phi^{\cal K}(S)$ is a good subset.

\vskip.2cm
     \noi      {\bf Lemma 9}  of 
Simon (1999): If $A$ is a good  subset of a cell $C$ 
              and if   $A\cap F$ is closed for every  $F\in {\cal
          P}^j$ with  
            $A\cap F\not= \emptyset$, then 
          $A=C$. \vskip.2cm

\noi 
We need a few more facts about $\Omega (X,J)$ for non-empty $X$ and $J$. If 
   $|J|\geq 2$ then  $\Omega (X,J)$  is
 topologically equivalent to a Cantor set, (Fagin, Halpern and Vardi
 1991).
  A Cantor set with the usual
            topology is a metric space. 
 Second  we can 
 perceive a Cantor set  as $\{ 0,1\} ^{\bf N}$, where 
 each finite sequence $a=(a^1, a^2, \dots , a^n)$  defines 
 a   cylinder subset $C(a)$ of $\{ 0,1\} ^{\bf N}$ by 
$C(a):= \{ z \in \{ 0,1\}^{\bf N}\
 | \ z^k =a^k \ \forall   k\leq n\}$. Furthermore 
 all  cylinder subsets are themselves topologically equivalent 
 to Cantor sets, and the same holds for finite unions of cylinder sets.
  Third, if $|J|\geq 2$ then   
 there exists  an uncentered 
cell of $\Omega (X,J)$ of finite fanout 
 that is dense in $\Omega (X,J)$  (Simon 1999).

\vskip.2cm 
\noi 
 Due to   topological formulations 
 of the centered property, to demonstrate that there is a 
 surjective cell without finite fanout  
   requires some topological insight.  Central to our 
 proof is 
     Theorem 9 of Chapter 
 12 of  E. Moise, (1977):  

 \vskip.2cm 
\noi  Let $X$ and $Y$ be two totally disconnected,  
  perfect, compact metric spaces (equivalently Cantor sets)  
  and let $X'$ and $Y'$ be countable    and 
   dense subsets of $X$ and $Y$, respectively.  There is a homeomorphism 
    between $X$ and $Y$ that is also a bijection  between  $X'$ and $Y'$. 
     \vskip.2cm

\noi 
We call a partition ${\cal P}$ of a metric  space $D$ 
         upper (respectively lower) hemi-continuous if the set valued  
         correspondence that maps every $d\in D$ to the partition 
          member of ${\cal P}$ containing $d$ is an upper (respectively 
           lower) hemi-continuous correspondence.   (We follow the definitions 
            of Klein and Thompson, 1984.)

 \vskip.2cm
 \noi  {\bf Lemma 1:} If          
       ${\cal K}:= (S;J; 
       ({\cal P}^j\ | \ j\in J); X;\psi)$  
   is a Kripke structure  
 with a topology (not necessarily that induced by the 
 formulas) such that\newline     
    1) for every $z\in \{ 0,1\} ^X$ the set 
 $\psi ^{-1}  (z)$ is  clopen (closed and open)   and\newline 2)  
   for every $j\in N$ 
   the partition ${\cal P}^j$ is lower and upper hemi-continuous, \newline 
     then the map $\phi^{\cal K} :S \rightarrow  
     \Omega (X,J)$ is continuous. 

            \vskip.2cm
  \noi {\bf Proof:}
It suffices to show that $\alpha ^{\cal K} (f)$ is 
 a clopen set for every $f\in {\cal L}(X,J)$.  
 We proceed by induction on the 
  structure of formulas. The
 claim is true for all $f=x\in X$ by  
   hypothesis. Due to the clopen property being closed 
     under complementation and finite intersection,
 it is likewise true for $\neg f$ and $f\wedge g$ if 
    it is true for $f$ and $g$.  For 
     some $f\in {\cal L}(X,J)$ 
      we assume that 
      $\alpha ^{\cal K} (f)$
is a clopen set.  $\alpha ^{\cal K} (k_j f)$ 
 is an open set by the upper semi-continuity of ${\cal P}^j$ and 
  the openness of $\alpha ^{\cal K} (f)$. 
  $S\backslash \alpha ^{\cal K} (k_j f)=
  \alpha ^{\cal K} (\neg k_j f)$ is an open 
   set by the openness of 
   $S\backslash \alpha ^{\cal K} ( f)$ 
    and the lower semi-continuity of ${\cal P}^j$. 
    \hfill $\Box$
                       \vskip.2cm

\noi {\bf Lemma 2:} Given $X$ and $J$ finite, for every $j\in J$ the partition 
 ${\cal Q}^j(X,J)$ of  $\Omega (X,J) $ is upper and lower hemi-continuous with
 respect to the topology induced by the formulas. 
\vskip.2cm 
\noi 
{\bf Proof:} 
Let $z_1, z_2, \dots $ be a sequence  of points in $\Omega (X,J)$ 
  converging to some  $z\in
P\in {\cal Q}^j$ with  
 $z_i\in P_i \in {\cal Q}^j$ for every $i=1,2,\dots$. 
\vskip.2cm 

\noi 
 To prove that
 ${\cal Q}^j$ is upper hemi-continuous it suffices  to show that if 
 $y_1, y_2, \dots $ is a sequence of   points in $\Omega (X,J)$ 
  converging to $y$ 
 with $y_1\in P_1, y_2 \in P_2, \dots$  
  then 
 $y$ is in $ P$. Let $f$ be any formula such that $k_j f\in y$. Since
 the $y_i$ converge to $y$  there is an $N$ such that
 for every $i\geq N$ it must hold that $k_j f$ is in both $y_i$ and 
 $z_i$. But this
 means that $k_j f$ is also in $z$. The same argument holds for 
  the formula $\neg k_j f$.   
\vskip.2cm 

\noi 
To prove that ${\cal Q}^j$ is lower hemi-continuous it suffices  to show
that if 
 $y\in P$ then there is a sequence of $y_1, y_2, \dots$ in $P_1, P_2,
 \dots$ respectively that converges to $y$.
 Because there are only countably many formulas and one can create 
 a new sequence from the diagonal of  sequences which come closer and
 closer to  $y$,   
  if the claim were not true then  there would be 
 some formula $f$ in  $y$ and an $N$  
 such that $f$ is not in any member of $P_i$ for all $i\geq N$. 
This would imply also that $k_j ( \neg f)$ is in $z_i$ for all $i\geq
N$ and likewise that $k_j(\neg f) $ is in $z$. 
 But this would contradict that the assumption that $f$ is in  $y$
 and $y$  is in $ P$. 
\hfill $\Box$.

\section{The Example}

Let $S$ equal $\Omega(X,\{ 1,2\})$ with  $X$ any finite non-empty set. 
 Let $C$ be an uncentered 
 cell of finite fanout that is dense in $S$. We assume that 
  $\pi: S \rightarrow  \{ 0,1\} ^{\bf N}$
 is a homeomorphism.  For every $n\in {\bf N}$ define 
$\pi_n: S \rightarrow \{ 0,1\} ^n$ by 
$\pi _n (z)$ equalling the 
   $a=(a^1, a^2, \dots , 
 a^n)\in \{ 0,1\} ^{n}$ such that   
 $\pi (z)=( a_1, \dots , a_n , \dots)$. This means that  
  $\pi_n^{-1}  \circ
 \pi_n (z) $ equals $C(\pi_n (z))$, the corresponding 
cylinder set. If $a$ is the empty sequence in 
 $\{ 0,1\} ^0$ then define $\pi_0 (z):= a$ and $\pi ^{-1} _0 \circ
 \pi_0 (z) = S$ for all $z\in S$. 
\vskip.2cm 

\noi 
 Let $z$ be any member 
 of $C$ and for every $i=1,2,\dots$ let $z_i$ be a member of $C$ 
 such that $\pi_{2i-2}(z_i)=\pi_{2i-2}(z)$   
  but $\pi_{2i}(z_i)\not= \pi_{2i}(z)$. 
  For every $i=1,2,\dots $ define non-empty and mutually disjoint sets 
 $A_{i,1}, A_{i,2}, \dots A_{i,i}$ in the following way.  
Let $A_{1,1}$ equal $ S \backslash \ (\pi^{-1}_2 \circ \pi_2 (z_1)\cup 
\pi^{-1}_2  \circ \pi_2 (z))$. For $1\leq k<i$ let 
 $A_{i,k}:=\pi_{2i-2}^{-1} \circ \pi_{2i-2}(z_k)\backslash \ 
\pi_{2i}^{-1}\circ \pi_{2i}(z_k)$ and 
 let $A_{i,i}:= \pi^{-1}_{2i-2}  \circ \pi_{2i-2} (z)\backslash \  
(\pi^{-1}_{2i} \circ \pi_{2i} (z_i)\cup \pi^{-1}_{2i}  \circ \pi_{2i} (z))$.  
 Because for every $a\in \{ 0,1\}^{2i}$ there are four members $b$ of 
 $\{ 0,1\}^{2i+2}$ such that $a=\pi_{2i} \circ \pi^{-1}_{2i+2}(b)$, 
  all the 
 sets $ A_{i,j}$ are non-empty and
  homeomorphic  to Cantor sets.  By Proposition 1, for every $i\geq 1$ and 
 $1\leq k\leq i$ there is a homeomorphism 
 $f_k: A_{i,1}\rightarrow A_{i,k}$ such that 
 $f_k$ maps $C\cap A_{i,1}$ bijectively to $C\cap A_{i,k}$. This implies for 
 every $i\geq 1$ that there exists an upper and 
 lower semi-continuous partition ${\cal P}^i$ of $C\cap (\cup_{k=1}^i 
A_{i,k})$ such that every partition member of ${\cal P}^i$ has 
 $i$ members, one member in $A_{i,k}$ for every $1\leq k\leq i$. 
 Notice that all the $A_{i,k}$ are mutually disjoint, meaning 
that $A_{i,k}=A_{i^*, j^*}$ if and only if $i=i^*$ and $k=k^*$.  Furthermore 
 the disjoint union $\cup_{i\geq 1} \cup _{1\leq k\leq i}A_{i,k}$ is 
 equal to $S \backslash \ \{ z, z_1, z_2, \dots\}$. 
Let ${\cal P}$ be $(\cup_{i=1}^{\infty} {\cal P}^i) 
\cup \{z,z_1, z_2, \dots\}$, 
 a partition of $C$. It is straightforward to check that ${\cal P}$ is 
 upper and lower semi-continuous.  We define 
 $\cal A$ be the Kripke structure   
 $(C;\{ 1,2, 3\};{\cal Q}^1|_C, {\cal Q}^2|_C,
{\cal P}; X,  \psi|_C)$, with the partition ${\cal P}$ corresponding
 to the third person.
\vskip.2cm 

\noi 
{\bf Theorem:}  $\phi^{\cal A}$ maps $C$ bijectively to 
 a cell of $\Omega (\{ 1,2,3\})$ that is surjective but without 
 finite fanout. \vskip.2cm
\noi 
{\bf Proof:}     
 We have by Lemma 1 that 
 $\phi^{{\cal A}}:C\rightarrow \Omega(X,\{ 1,2, 3\})$ is continuous. 
 Since every member of 
${\cal Q}^1|_C$, ${\cal Q}^2|_C$,
or ${\cal P}$ is compact, their images in $\Omega (X,\{ 1,2,3\})$ 
are also 
 compact.    By Lemma 9 of Simon (1999) 
   $\phi^{\cal A}$ maps
 $C$ surjectively to a 
 cell $\phi^{\cal A}(C)$ of $\Omega (X,\{ 1,2,3\})$. 
        Between any two points of
 $\phi^{\cal A}(C)$ there is an adjacency      
         path using images of members of 
         ${\cal Q}^1|_C$ and ${\cal Q}^2|_C$,  therefore  there can be 
                  no proper good subset of 
   $\phi^{\cal A}(C)$. By Lemma 7 of Simon (1999) 
 this implies that $\phi^{\cal A}(C)$ is 
 a surjective cell.
             Since for every 
             $f\in {\cal L}(X,\{ 1, 2\})$ 
              $\alpha ^{\Omega (X,\{ 1,2\})}(f)$       
              gets mapped to $\alpha ^{\Omega (X,\{ 1,2,3\})}(f)$, 
               $\phi^{\cal A}$ is an injective and an open map
 (meaning that open sets are mapped to open sets), and 
 therefore  the map $\phi^{\cal A}$ is also  
 a homeomorphism of $C$ to  $\phi^{\cal A}(C)$. 
               Therefore the
 image of the one infinite  set in ${\cal P} $ is 
                also an infinite 
 set in the cell $\phi^{\cal A}(C)$, 
 which implies that this cell of 
 $\Omega (X,\{ 1,2,3\})$ does not have finite fanout.
 \hfill q.e.d. \vskip.2cm

  \section {Acknowledgements} 
       Many helped in finding the best citation for  Proposition 1; the  
       lemma  has many interesting variations, including the same conclusion 
        for open intervals proven  by G. Cantor  (1895).

\section {References}
\begin{description}
 \item  [Aumann, R. (1976),] ``Agreeing to Disagree," Annals of Statistics 4, 
 pp. 1236-1239. 
 \item [Bacharach, M., Gerard-Varet, L.A., Mongin, P, and Shin, H., eds.] 
  {\bf (1997),} {\em Epistemic Logic and the 
  Theory of Games and Decisions}, Dordrecht, 
   Kluwer.
 \item [Cantor, G. (1895),] ``Beitraege zur Begruendung der transfiniten 
  Mengenlehre," Mathematische Annalen, Band 46, pp. 481-512.
 \item [Fagin, R. (1994),] ``A Quantitative Analysis of Modal Logic," Journal 
  of Symbolic Logic 59, pp. 209-252.
 \item [Fagin, R., Halpern, Y.J. 
and Vardi, M.Y. (1991),] ``A Model-Theoretic 
 Analysis of Knowledge," Journal of the A.C.M. 91 (2), pp. 382-428.
    \item[Halpern, J., Moses, Y. (1992),] ``A Guide to Completeness and 
     Complexity for Modal Logics of Knowledge and Belief," 
     Artificial Intelligence 54, pp. 319-379.
 \item [Heifetz, A. and Samet, D. (1998),] ``Knowledge Spaces with Arbitrarily 
     High Rank," 
   Games and Economic Behavior 22, No.2, pp. 260-273. 
   \item [Heifetz, A. and Samet, D. (1999),] ``Hierarchies of Knowledge: 
 An Unbounded Stairway", Mathematical Social Sciences, Vol. 38, pp. 157-170.
   \item [Hughes, G.E., Cresswell, M.J. (1968),] 
 {\em An Introduction to Modal Logic,} Routledge.
 \item  [Klein, E., Thompson, A.  
(1984),]  {\em Theory of Correspondences,} Wiley. 
\item  [Lewis, D.  
(1969),]  {\em Convention: A Philosophical Study,} Harvard University Press. 
 \item [Moise, E. (1977),] {\em Geometric Topology 
 in Dimensions Two and Three,} 
 Graduate Texts in Mathematics 47, Springer.
 \item [Rubinstein, A. (1989)] `` The Electronic Mail Game:
 A Game with Almost Common Knowledge'', American Economic
 Review, Vol. 79, pp. 385-391.
 \item [Samet, D. (1990)] ``Ignoring Ignorance and Agreeing to Disagree," 
  Journal of Economic Theory, Vol. 52, No. 1, pp. 190-207.
\item [Simon, R. (1999)] 
``The Difference between Common Knowledge of Formulas and Sets", International
 Journal of Game Theory, Vol. 28, No. 3, pp. 367- 384.
\item [Simon, R. (2003)] ``Games of Incomplete Information, Ergodic Theory, 
 and the Measurability of Equilibria'', Israel Journal of Mathematics, 
 Vol. 138, pp. 73-92.

 \end {description}

\end{document}